\documentclass[12pt]{article}

\usepackage{amsmath,amssymb}
\usepackage[dvips]{graphicx}
\usepackage{subfigure}
\usepackage{color}
\usepackage{dcolumn}
\usepackage{bm}
\usepackage[numbers,super,comma,sort&compress]{natbib}
\usepackage{multirow}
\usepackage{fancyhdr}
\usepackage{latexsym}
\usepackage{amscd}
\usepackage{cases}
\usepackage{algorithmic}
\usepackage{algorithm}
\usepackage{booktabs}
\usepackage{enumerate}

\definecolor{background-color}{gray}{0.98}

\renewcommand\AA{\text{\r A}}

\usepackage{color}

\begin{document}

\title{Efficient and Qualified Mesh Generation for Gaussian Molecular Surface Using Piecewise Trilinear Polynomial Approximation
}
\author{
Tiantian Liu\thanks{State Key Laboratory of Scientific  and Engineering Computing, Institute of Computational Mathematics and Scientific  Engineering Computing, Academy of Mathematics and Systems Science, Chinese Academy of Sciences, Beijing 100190, China. Corresponding author: Benzhuo Lu (E-mail: bzlu@lsec.cc.ac.cn).
},
Minxin Chen\thanks{Center for System Biology, Department of Mathematics, Soochow University, Suzhou 215006, China. Corresponding author: Minxin Chen (E-mail: chenminxin@suda.edu.cn).},
Benzhuo Lu{$^*$,}
}

\date{}

\maketitle
\begin{abstract}
Recent developments for mathematical modeling and numerical simulation of biomolecular systems raise new demands for qualified, stable, and efficient surface meshing, especially in implicit-solvent modeling\cite{Lu08}. In our former work, we have developed an algorithm for manifold triangular meshing for large Gaussian molecular surfaces, TMSmesh\cite{TMS11,TMS12}. In this paper, we present new algorithms to greatly improve the meshing efficiency and qualities, and implement into a new program version, TMSmesh 2.0.
In TMSmesh 2.0, in the first step, a new adaptive partition and estimation algorithm is proposed to locate the cubes in which the surface are approximated by piecewise trilinear surface with controllable precision.
Then, the piecewise trilinear surface is divided into single valued pieces  by tracing along the fold curves, which ensures that the generated surface meshes are manifolds.
Numerical test results show that TMSmesh 2.0 is capable of handling arbitrary sizes of molecules  and achieves ten to hundreds of times speedup over the previous algorithm. The result surface meshes are manifolds and can be directly used in boundary element method (BEM) and finite element method (FEM) simulation.
The binary version of TMSmesh 2.0 is downloadable at the web page http://lsec.cc.ac.cn/$\sim$lubz/Meshing.html.
\end{abstract}

\section{Introduction}
Molecular surface mesh generation is a prerequisite for using boundary element method (BEM) and finite element method (FEM) in the implicit-solvent modeling (e.g., see a review in \cite{Lu08}).
Recent developments in implicit-solvent modeling of biomolecular systems raise new demands for qualified, stable, and
efficient surface meshing. Main concerns for improvement on existing methods for molecular surface mesh generation are efficiency, robustness, and mesh quality.
Efficiency is necessary for simulations/computations requiring frequent mesh generation or requiring meshing for large systems. Robustness here means the meshing method is stable and can treat various, even arbitrary, sizes of molecular systems within computer power limitations. Mesh quality
relates to mesh smoothness (avoiding sharp solid angles, etc.), uniformness (avoiding elements with very sharp angles or zero area), topological correctness (manifoldness, avoiding isolated vertices, element intersection, single-element-connected edges,
etc.) and fidelity (faithful to the original defined molecular surface). The quality requirement is critical for some numerical techniques, such as finite element method, to achieve converged and reasonable results, which makes it a more demanding task in this aspect than the mesh generations only for the purposes of visualization or some structural geometry analysis.

There are various kinds of definitions for molecular surface, including the van der Waals (VDW) surface, the solvent accessible surface (SAS)\cite{Lee71}, the solvent excluded surface (SES)\cite{Richards77}, the minimal molecular surface\cite{Bates08}, the molecular skin surface\cite{Edelsbrunner99} and the Gaussian surface.
The VDW surface is defined as the surface of the union of the spherical atomic surfaces with VDW radius of each atom within the molecule.
The SAS and SES are represented by the trajectory of the center and the inter-boundary of a rolling probe on the VDW surface, respectively.
The minimal molecular surface is defined as the result of the minimization of a type of surface energy.
The molecular skin surface is the envelope of an infinite family of spheres derived from atoms by convex combination and shrinking.
The Gaussian surface is defined as a level set of the summation of Gaussian kernel functions:
\begin{equation}
\{ \vec x \in R^3,\phi \left(\vec x \right) = c \},
\end{equation}
where
\begin{equation}
\phi \left(\vec x \right) = \sum _{i=1}^{N}e^{-D(\Vert \vec x - \vec x_i \Vert ^2 -r_i^2)},
\label{Gaussian1}
\end{equation}
$\vec x_i $ and $r_i $ are the location and radius of the $i$th atom.
$D$ is the decay rate of the Gaussian kernel.
$c$ is the isovalue and it controls the volume enclosed by the Gaussian surface.
These two parameters, $D$ and $c$ can be chosen properly to make the Gaussian surface approximate the SES, SAS and VDW surface well.\cite{liutt15}

For SAS and SES, numerous works have been committed to the computation of the molecular surface in the literature.
In 1983, Connolly proposed algorithms to calculate the molecular surface and SAS analytically.\cite{Connolly83a,Connolly83b}
In 1995, a popular program, GRASP, for visualizing molecular surfaces was presented.\cite{Nicholls95}
An algorithm named SMART  for triangulating SAS into curvilinear elements was proposed by Zauhar\cite{Smart95}.
The software MSMS was proposed by Sanner et al. in 1996 to mesh the SES  and  is a widely used program for molecular surface triangulation due to its high efficiency.\cite{ Sanner96}
In 1997, Vorobjev et al. proposed SIMS, a method of calculating a smooth invariant molecular dot surface, in which an exact method for removing self-intersecting parts and smoothing the singular regions of the SES was presented.\cite{Vorobjev97b}
Ryu et al. proposed a method based on beta-shapes\cite{Ryu07}, which is a generalization of alpha shapes\cite{Edelsbrunner94}.
Can et al. proposed LSMS to generate the SES on grid points using level-set methods.\cite{Can06}
In 2009, a program, EDTsurf, based on LSMS was proposed used for generating the VDW surface, SES and SAS.\cite{EDT09}
A ray-casting-based algorithm, NanoShaper, is proposed to generate SES, skin surface and Gaussian surface in 2013.\cite{Decherchi13_2}

For skin surface, Chavent et al. presented MetaMtal to visualize the molecular skin surface using ray-casting method \cite{Chavent2008209}, and Cheng et al. used restricted union of balls to generate mesh for molecular skin surface\cite{Cheng2009196}.
For minimal surface, Wei et al. \cite{Bates08} constructed a surface-based energy functional, and use minimization and isosurface extraction processes to obtain a so-called minimal molecular surface.

For the Gaussian surface, existing techniques for triangulating an implicit surface can be used to mesh the Gaussian surface.
These methods are divided into two main categories: spatial partition and continuation methods.
The well known marching cube method \cite{Lorensen87} and dual contouring method \cite{Ju02} are examples of the spatial partition methods.
In 2006, Zhang et al. \cite{Zhang06} used a modified dual contouring method to generate meshes for biomolecular structures.
A later tool, GAMer \cite{Yu08b}, was developed for both the generation and improvement of the Gaussian surface meshes.
An efficient mesh generation algorithm accelerated by multi-core CPU and GPU was also proposed in 2013.\cite{Zhang13}

Most of those software have some issues according to above mentioned criteria for mesh generation, e.g., MSMS and GAMer generate many non-manifold defects in the mesh, fidelity is not well preserved for the EDTsurf and GAMer surfaces, and for the marching cube and grid-based methods the memory requirements can be huge when treating large molecules. More detailed comparison and discussion of the software can be found in \cite{liutt15}.
As MSMS is a most commonly used software in this area, we will still use it as a main reference for our new algorithm in this article.
In 2011, we have proposed an algorithm and implemented in the program TMSmesh for triangular meshing of the Gaussian surface.\cite{TMS11,TMS12,liutt15}
The trace technique which is a generalization of adaptive predictor-corrector technique is used in TMSmesh to connect sampled surface points.
TMSmesh contains two steps.
The first step is to compute the intersection points between the molecular Gaussian surface and the lines parallel to $x$-axis.
In the second step, the sampled surface points are connected through three algorithms to form loops, and the whole closed manifold surface is decomposed into a collection of patches enclosed by loops on the surface.
The patches are finally single valued on $x, y, z$ directions, so these pieces can be treated as 2-dimensional polygons and be easily triangulated through standard triangulation algorithms.
In TMSmesh, there are no problems of overlapping, gap filling, and selecting seeds that need to be considered in traditional continuation methods.
TMSmesh performs well in the following aspects.
Firstly, TMSmesh is robust.
TMSmesh succeeds to generate surface meshes for biomolecules comprised of more than one million atoms.
Secondly, the meshes produced by TMSmesh have good qualities (uniformness, manifoldness).
Thirdly, the generated surface mesh preserves the original molecular surface features and properties (topology, surface area and enclosed volume, and local curvature).
However, as to the aspect of computational efficiency, although the computational complexity is linear with respect to the number of atoms as shown in \cite{TMS12}, the overall low efficiency of TMSmesh still needs to be improved.

In this paper, we proposed a new algorithm and updated program version TMSmesh 2.0 to mesh the Gaussian surface efficiently.
Firstly, the space are adaptively divided into cubes and an algorithm of dividing cubes and estimating the error between the Gaussian surface and approximated  trilinear polynomial in each cube is developed. With this algorithm, the Gaussian surface is  approximated by piecewise trilinear surface. Then, in each cube, the trilinear surface is divided into a collection of single valued pieces on $x, y, z$ directions by tracing along fold curves (in the trilinear surface case the fold curve can be directly calculated analytically). Finally, each single valued piece is triangulated by ear clipping algorithm\cite{98triangulation,meisters1975polygons}.

This paper is organized as follows.
The new algorithm for triangulating the Gaussian surface is introduced in the Meshing Algorithm Section.
In the Experimental Results Section, some examples and applications are presented.
The final section, Conclusion, gives some concluding remarks.

\section{Meshing Algorithm}\label{sec1}
\subsection{Algorithm Outline}
In this section, we describe the algorithms to construct the triangular surface meshes.
The inputs of our method are PQR files which contains a list of centers and radii of atoms.
The output of our method are OFF files which contains the triangular meshes.
Our algorithm contains two stages, the first stage is an adaptive estimation and division process. The Gaussian surface is approximated by piecewise trilinear surface within controllable error.
The second stage is to partition each piece of trilinear surface into single valued patches along $x, y, z$ directions by tracing along the fold curves. Then each single valued patch is triangulated by the ear clipping algorithm\cite{98triangulation,meisters1975polygons}.
In the following subsections, each stage is described in detail.

\subsection{Approximating the Gaussian Surface by Piecewise Trilinear Surface}

In this stage, the space is divided into cubes adaptively and in each final cube, the Gaussian surface is close to a trilinear surface whose error is controllable.
Initially, the molecule is placed in a three-dimensional orthogonal grid consisting $n_x \times n_y \times n_z$ cubes.
The initial grid is very coarse.
Then the grid is refined adaptively by the following estimation and division steps.
\begin{itemize}
  \item Step 1, In each cube, $\phi(x,y,z)$ is approximated by a \emph{n}th-degree polynomial $\tilde P(x,y,z)$, i.e., the Gaussian surface $\phi(x,y,z) = c$ is replaced by the polynomial surface $\tilde P(x,y,z) = c$.
  \item Step 2, the lower and the upper bound of $\tilde P(x,y,z)$, denoted by $L$ and $U$, in each cube is estimated. If the isovalue $c$ belongs to $[L,U]$, the cube has intersection with the surface $\tilde P(x,y,z) = c$ and we go to Step 3, otherwise, the cube is abandoned.
  \item Step 3, divide each left cube into 8 smaller child cubes, and compute the expression of $\tilde P(x,y,z)$ in each child cube. When the child cubes become smaller, the coefficients of higher order terms (higher than the linear order) of $\tilde P(x,y,z)$ go to zero. If they are under some user-specified bound, approximate $\tilde P(x,y,z)$ by trilinear polynomial, otherwise, go to Step 2.
\end{itemize}

With above processes, the Gaussian surface finally is approximated by piecewise trilinear surfaces in cubes with different sizes. In the following subsections, we explain the details of above estimation and division process.
\subsubsection{Approximation with \emph{n}th-degree polynomial}
Firstly, without loss of generality, we only consider the case of $D=1$, then eq \eqref{Gaussian1} is written into the following one:
\begin{equation}
\phi \left(\vec x \right) = \sum _{i=1}^{N}e^{-(\Vert \vec x - \vec x_i \Vert ^2 - r_i^2 )}=\sum _{i=1}^{N}e^{r_i^2}e^{-\left( x-x_i \right)^2}e^{-\left( y-y_i \right)^2}e^{-\left( z-z_i \right)^2}.
\label{Gaussian2}
\end{equation}
In an arbitrary cube $[a,b]\times [c,d]\times [e,f]$, eq \eqref{Gaussian2} can be approximated by
\begin{equation}
P(x,y,z) = \sum _{i=1}^{N}e^{r_i^2}P_n(x,x_i,a,b)Q_n(y,y_i,c,d)R_n(z,z_i,e,f),
\label{poly1}
\end{equation}
where
\begin{equation}
P_n(x,x_i,a,b) = \sum_{j=0}^n \alpha_j(x_i,a,b)L_j(\frac{2x-(a+b)}{b-a}),
\label{poly2_1}
\end{equation}
\begin{equation}
Q_n(y,y_i,a,b) = \sum_{j=0}^n \beta_j(y_i,c,d)L_j(\frac{2y-(c+d)}{d-c}),
\label{poly2_2}
\end{equation}
\begin{equation}
R_n(z,z_i,e,f) = \sum_{j=0}^n \gamma_j(z_i,e,f)L_j(\frac{2z-(e+f)}{f-e}),
\label{poly2_3}
\end{equation}
and
$$\alpha_j = \frac{2}{b-a}\int_{a}^b \phi(\frac{2x-(a+b)}{b-a})L_j(\frac{2x-(a+b)}{b-a})\, dx,$$
$$\beta_j = \frac{2}{d-c}\int_{c}^d \phi(\frac{2y-(c+d)}{d-c})L_j(\frac{2y-(c+d)}{d-c})\, dy,$$
$$\gamma_j = \frac{2}{f-e}\int_{e}^f \phi(\frac{2z-(e+f)}{f-e})L_j(\frac{2z-(e+f)}{f-e})\, dz,$$
$L_j(\cdot)$ is Legendre polynomial of order $j$ and \emph{n} is set as 3 in our work.
However, $P(x, y, z)$ is not continuous between neighbored cubes, so we do the following corrections of $P(x, y, z)$ to make $P(x, y, z)$ be $C^0$ continuous in the whole domain.
For one component $P_n(x,x_i,a,b)$, we introduce two variables $\epsilon_0(x_i,a,b) $ and $\epsilon_1(x_i,a,b)$ as follows.
\begin{equation}
\begin{aligned}
\tilde{P}_n(x,x_i,a,b) &= P_n(x,x_i,a,b) + \epsilon_0(x_i,a,b)L_{n-1}(\frac{2x-(a+b)}{b-a})+ \epsilon_1(x_i,a,b)L_n(\frac{2x-(a+b)}{b-a}) \\
 &= \alpha_0(x_i,a,b)L_0(\frac{2x-(a+b)}{b-a})+\cdots +\alpha_n(x_i,a,b)L_n(\frac{2x-(a+b)}{b-a}) \\
&+\epsilon_0(x_i,a,b)L_{n-1}(\frac{2x-(a+b)}{b-a})+ \epsilon_1(x_i,a,b)L_n(\frac{2x-(a+b)}{b-a}).
\end{aligned}
\label{poly3}
\end{equation}
The following two equations make $\tilde P_n(x,x_i,a,b)$ equal $x$ component of $\phi\left(\vec x \right)$ on the boundary of the box and be $C^0$ continuous along $x$ directions.
\begin{subequations}
\begin{numcases}{}
\tilde P_n(a,x_i,a,b) = e^{-(a-x_i)^2}\\
\tilde P_n(b,x_i,a,b) = e^{-(b-x_i)^2}.
\end{numcases}
\label{condition1}
\end{subequations}
$\epsilon_0(x_i,a,b) $ and $\epsilon_1(x_i,a,b) $ can be easily solved from eq (\ref{condition1}).
Then $\tilde P_n(x,x_i,a,b)$ is written as the following one:
\begin{equation}
\tilde P_n(x,x_i,a,b) = \sum_{j=0}^n \tilde \alpha_j(x_i,a,b)L_j(\frac{2x-(a+b)}{b-a}),
\label{pn}
\end{equation}
where
\begin{equation}\label{alpha}
  \tilde \alpha_j(x_i,a,b) =\begin{cases}
  \alpha_j(x_i,a,b) & \text{$j < n-1$}, \\
  \alpha_{n-1}(x_i,a,b) + \epsilon_0(x_i,a,b) & \text{$j=n-1$}, \\
   \alpha_n(x_i,a,b) + \epsilon_1(x_i,a,b) & \text{$j=n$}.
   \end{cases}
\end{equation}

After above correction, $\tilde P_n(x,x_i,a,b)$ is the best least square approximation for $x$ component of $\phi(\vec x)$ in the space spanned by \{$L_i$, $i=0,...,n-2$\} and it is also $C^0$ continuous on the boundaries of the cubes along $x$ direction.
The same method should be used to correct $Q_n(y,y_i,c,d)$ and $R_n(z,z_i,e,f)$ to make $\tilde P(x,y,z)$ be $C^0$ continuous along $y, z$ directions.
We have
\begin{equation}
 \tilde Q_n(y,y_i,a,b) = \sum_{j=0}^n \tilde \beta_j(y_i,c,d)L_j(\frac{2y-(c+d)}{d-c}),
\end{equation}
\begin{equation}
  \tilde R_n(z,z_i,e,f) = \sum_{j=0}^n \tilde \gamma_j(z_i,e,f)L_j(\frac{2z-(e+f)}{f-e}),
\end{equation}
where the forms of $\tilde \beta_j(y_i,c,d)$ and $\tilde \gamma_j(z_i,e,f)$ are similar to $\tilde \alpha_j(x_i,a,b)$ in eq \eqref{alpha}.
Then the new $n$-th polynomial is written as
\begin{equation}
\tilde P(x,y,z) = \sum _{i=1}^{N}e^{r_i^2}\tilde P_n(x,x_i,a,b)\tilde Q_n(y,y_i,c,d)\tilde R_n(z,z_i,e,f)
\label{polynew}
\end{equation}
for $x \in[a,b], y \in [c,d], z\in[e,f]$. In practical computation of $\tilde P(x,y,z)$, we only need to compute the summation in eq \eqref{polynew} with respect to the  neighborhood $\{x_i,y_i,z_i\}$ of the cube $[a,b]*[c,d]*[e,f]$, since the  kernel $e^{-||\vec{x}-\vec{x}_i||^2}$ decay very quickly when $||\vec{x}-\vec{x}_i||$ goes to large.

\subsubsection{Estimation of upper and lower bound of $\tilde P(x,y,z)$}\label{estimate}
In order to rule out the cubes having no surface points, the lower and upper bound of $\tilde P(x,y,z)$ in the cube is estimated.
$\tilde P(x,y,z)$ in eq \eqref{polynew} can be written in the form of product of tensor:
\begin{equation}
\tilde P(x,y,z) = A \bar \times_1 \vec L(\frac{2x-(a+b)}{b-a}) \bar \times_2 \vec L(\frac{2y-(c+d)}{d-c})\bar \times_3 \vec L(\frac{2z-(e+f)}{f-e}),
\label{tensor1}
\end{equation}
where
$A = \sum _{i=1}^{N}e^{r_i^2}B_i$ and
\begin{subequations}
\begin{numcases}{}
B_i=b_i^{(1)}\otimes b_i^{(2)} \otimes b_i^{(3)}\\
b_i^{(1)}=(\tilde \alpha_0(x,x_i,a,b),\tilde \alpha_1(x,x_i,a,b),\cdots ,\tilde \alpha_n(x,x_i,a,b))\\
b_i^{(2)}=(\tilde \beta_0(y,y_i,a,b),\tilde \beta_1(y,y_i,a,b),\cdots ,\tilde \beta_n(y,y_i,a,b))\\
b_i^{(3)}=(\tilde \gamma_0(z,z_i,a,b),\tilde \gamma_1(z,z_i,a,b),\cdots ,\tilde \gamma_n(z,z_i,a,b))\\
\vec L(\frac{2x-(a+b)}{b-a})=(L_0(\frac{2x-(a+b)}{b-a}),\cdots ,L_n(\frac{2x-(a+b)}{b-a}))\\
\vec L(\frac{2y-(c+d)}{d-c})=(L_0(\frac{2y-(c+d)}{d-c}),\cdots ,L_n(\frac{2y-(c+d)}{d-c}))\\
\vec L(\frac{2z-(e+f)}{f-e})=(L_0(\frac{2z-(e+f)}{f-e}),\cdots ,L_n(\frac{2z-(e+f)}{f-e})).
\end{numcases}
\label{tensor2}
\end{subequations}

$\otimes$ is the product of tensor. $\bar \times_k$ is $k$-mode (vector) product of a tensor $X \in \mathbb{R}^{I_1\times I_2 \times I_3}$ with a vector $V \in \mathbb{R}^{I_k}$ denoted by $X \bar \times_k V$ and is of size $I_1\times \cdots \times I_{k-1} \times I_{k+1} \times \cdots \times I_3$. \cite{tensor} Its $i_1\cdots i_{k-1}i_{k+1}\cdots i_3$ entry is as follows.
\begin{equation}
(X \bar \times_k V)_{i_1\cdots i_{k-1}i_{k+1}\cdots i_3} = \sum_{i_k=1}^{I_k}x_{i_1}x_{i_2} x_{i_3}v_{i_k}.
\label{nproduct}
\end{equation}
$A$ is a three-dimensional tensor whose size is $(n+1)\times (n+1) \times(n+1)$, where $n$ is the degree of the polynomial $\tilde P(x,y,z)$.
To get the lower and upper bound of $\tilde P(x,y,z)$, firstly, the main part of $\tilde P(x,y,z)$ is obtained by doing singular value decomposition (SVD) for $A$.
Secondly, the upper and the lower bounds of the main part and the remainder are estimated respectively.

Here we use Singular Value Decomposition(SVD) for $A$ to approximate $\tilde P(x,y,z)$ by a multiplication of three polynomials in $x, y, z,$ respectively.
Taking $n = 3$ and $A=[a_{ij}]_{4 \times 4 \times 4}$ for example, the algorithm of SVD is as follows.

\textit{Step 1, transform A into a two-dimensional matrix.}
\addtocounter{MaxMatrixCols}{15}
\begin{equation}
A_1=
\begin{bmatrix}
a_{111} & \dots & a_{141} & a_{211} & \dots & a_{241} & a_{311} & \dots & a_{341} & a_{411} & \dots & a_{441} \\
a_{112} & \dots & a_{142} & a_{212} & \dots & a_{242} & a_{312} & \dots & a_{342} & a_{412} & \dots & a_{442} \\
a_{113} & \dots & a_{143} & a_{213} & \dots & a_{243} & a_{313} & \dots & a_{343} & a_{413} & \dots & a_{443} \\
a_{114} & \dots & a_{144} & a_{214} & \dots & a_{244} & a_{314} & \dots & a_{344} & a_{414} & \dots & a_{444}
\end{bmatrix}
\label{matrix1}
\end{equation}

\textit{Step 2, do singular value decomposition to $A_1$.}
\begin{equation}
A_1  = UDV^{\ast},
\label{SVD1}
\end{equation}
where $U = (\vec u_1, \vec u_2, \vec u_3, \vec u_4)$ is a $4 \times 4$ matrix, $V^\ast$ is a $4 \times 16$ matrix and
\begin{equation}
D=
\begin{bmatrix}
\sigma_1 & 0& 0&0 \\
0 & \sigma_2 &0 &0 \\
0 & 0 & \sigma_3 &0 \\
0 & 0 & 0 & \sigma_4
\end{bmatrix}
\label{SVD2}
\end{equation}
If $j$ satisfies
\begin{equation}
\text{min}\{j:\sum_{i=1}^{j}\sigma_i \geq 0.99\sum_{i=1}^{4}\sigma_i\},
\label{SVD3}
\end{equation}
we reserve $\sigma_1,\dots,\sigma_j$ and abandon $\sigma_{j+1}, \dots ,\sigma_4$.

\textit{Step 3, transform each row of $V^\ast$ into a square matrix.}
If
\begin{equation}
V^\ast=
\begin{bmatrix}
v_{1,1} & v_{1,2} & \dots & v_{1,16} \\
v_{2,1} & v_{2,2} & \dots & v_{2,16} \\
v_{3,1} & v_{3,2} & \dots & v_{3,16} \\
v_{4,1} & v_{4,2} & \dots & v_{4,16}
\end{bmatrix},
\label{SVD4}
\end{equation}
we can transform the $i$-th row of $V^\ast$ into a $4\times 4$ matrix denoted by $V_i$:
\begin{equation}
V_i=
\begin{bmatrix}
v_{i,1} & v_{i,2} & v_{i,3} & v_{i,4} \\
v_{i,5} & v_{i,6} & v_{i,7} & v_{i,8} \\
v_{i,9} & v_{i,10} & v_{i,11} & v_{i,12} \\
v_{i,13} & v_{i,14} & v_{i,15} & v_{i,16}
\end{bmatrix}.
\end{equation}

\textit{Step 4, do SVD for $V_i$ respectively.}

\begin{equation}
V_i = W_iD_iZ_i ,
\label{SVD5}
\end{equation}
where
$W_i=(\vec {w_1^i},\vec {w_2^i},\vec {w_3^i},\vec {w_4^i}),D_i=diag(d_1^i,d_2^i,d_3^i,d_4^i), Z_i=(\vec {z_1^i},\vec {z_2^i},\vec {z_3^i},\vec {z_4^i})^T$.
If $j_i$ satisfies $$\text{min}\{j_i:\sum_{k=1}^{j_i}d_k^i \geq 0.99\sum_{k=1}^{4}d_k^i\}, $$
$d_1^i,\dots,d_{j_i}^i$ are reserved and $d_{j_i+1}^i, \dots, d_4^i$ are abandoned.
Therefore, $V_i$ can be approximated by the following formula
\begin{equation}
V_i \approx d_1^i\vec{w_1^i} \otimes \vec{z_1^i} + \dots + d_{j_i}^i\vec{w_{j_1}^i} \otimes \vec{z_{j_i}^i} .
\label{SVD6}
\end{equation}

Through the above calculation, $A$ can be approximated by
\begin{equation}
A \approx \sum_{i=1}^{j}\sigma_i \vec u_i \otimes (\sum_{k = 1}^{j_i}d_i^k \vec{w_i^k} \otimes \vec{z_i^k}).
\label{SVD7}
\end{equation}
The summation at the right side of eq \eqref{SVD7} is denoted by $\tilde A$.
As a result, $A$ can be split by $A= \tilde A + R$, where $\tilde A$ is the main part and $R$ is the residue part.

With above SVD process, $\tilde P(x,y,z)$ can be converted to the following form:
\begin{equation}
\tilde P(x,y,z) = S(x,y,z)+T(x,y,z),\\
\label{poly5}
\end{equation}
where
\begin{equation}
S(x,y,z) = \tilde A \times_1 \vec L(\frac{2x-(a+b)}{b-a}) \times_2 \vec L(\frac{2y-(c+d)}{d-c}) \times_3 \vec L(\frac{2z-(e+f)}{f-e}),
\label{poly6}
\end{equation}
\begin{equation}
T(x,y,z) = R \times_1 \vec L(\frac{2x-(a+b)}{b-a}) \times_2 \vec L(\frac{2y-(c+d)}{d-c}) \times_3 \vec L(\frac{2z-(e+f)}{f-e}).
\label{poly7}
\end{equation}

For $S(x,y,z)$, the upper bound and lower bound are estimated through the following steps.
\begin{equation}
S(x,y,z) = \sum_{i=1}^{j}\sigma_i\tilde U^i(x)\left[\sum_{k=1}^{j_i}d_k^i \tilde W_k^i(y)\tilde Z_k^i(z) \right],
\label{Sxyz}
\end{equation}
where
\begin{subequations}
\begin{numcases}{}
\tilde U^i(x) = \vec u_i\cdot \vec L(\frac{2x-(a+b)}{b-a}) \\
\tilde W_k^i(y) = \vec {w_k^i} \cdot \vec L(\frac{2y-(c+d)}{d-c}) \\
\tilde Z_k^i(z) = \vec{z_k^i} \cdot \vec L(\frac{2z-(e+f)}{f-e}).
\end{numcases}
\end{subequations}

Firstly, we estimate the upper bound and lower bound of one dimensional polynomial $\tilde W_k^i(y)$ and $\tilde Z_k^i(z)$ respectively.
The upper and lower bound of $\tilde W_k^i(y)$ are denoted by $M_k^y$ and $m_k^y$.
And the upper and lower bound of $\tilde Z_k^i(z)$ are denoted by $M_k^z$ and $m_k^z$.
Secondly, the upper bound and lower bound of $\tilde W_k^i(y)\tilde Z_k^i(z)$ are estimated by
\begin{equation}
M_k^{yz} = max\{M_k^{y}M_k^{z}, M_k^{y}m_k^{z}, m_k^{y}M_k^{z}, m_k^{y}m_k^{z}\},
\end{equation}
\begin{equation}
m_k^{yz} = min\{M_k^{y}M_k^{z}, M_k^{y}m_k^{z}, m_k^{y}M_k^{z}, m_k^{y}m_k^{z}\}.
\end{equation}
Then the upper bound of $\sum_{k=1}^{j_i}d_k^i \tilde W_k^i(y)\tilde Z_k^i(z)$ is
\begin{equation}
M_i^{yz} = \sum_{k=1}^{j_i}d_k^iM_k^{yz}
\end{equation}
and the lower bound is
\begin{equation}
m_i^{yz} = \sum_{k=1}^{j_i}d_k^im_k^{yz}.
\end{equation}
Finally, we estimate the upper bound and the lower bound of $\tilde U^i(x)$ which are denoted by $M_i^x$ and $m_i^x$.
Then the bounds of $\tilde U^i(x)\left[\sum_{k=1}^{j_i}d_k^i \tilde W_k^i(y)\tilde Z_k^i(z) \right]$ can be estimated by
\begin{equation}
M_i = max\{M_i^{x}M_i^{yz},M_i^{x}m_i^{yz},m_i^{x}M_i^{yz},m_i^{x}m_i^{yz}\},
\end{equation}
\begin{equation}
m_i = min\{M_i^{x}M_i^{yz},M_i^{x}m_i^{yz},m_i^{x}M_i^{yz},m_i^{x}m_i^{yz}\}.
\end{equation}
Therefore, the upper bound and lower bound of $S(x,y,z)$ is
\begin{equation}
M = \sum_{i=1}^{j}\sigma_iM_i
\end{equation}
and
\begin{equation}
m = \sum_{i=1}^{j}\sigma_im_i.
\end{equation}

The range of each entry of $\vec L(\cdot)$ is $[-1,1]$.
Therefore, $T(x,y,z)$ can be estimated by
\begin{equation}
\begin{aligned}
|T(x,y,z)| & = |R \times_1 \vec L(\frac{2x-(a+b)}{b-a}) \times_2 \vec L(\frac{2y-(c+d)}{d-c}) \times_3 \vec L(\frac{2z-(e+f)}{f-e})| \\
& \leq |\sum R_{ijk}| \\
& \leq \sum |R_{ijk}|, \\
\end{aligned}
\label{bound1}
\end{equation}
where $R_{ijk}$ is the $(i,j,k)$ entry of $R$.

As a result, the upper and lower bound of $\tilde P(x,y,z)$ is
\begin{equation}
U=M+\sum |R_{ijk}|,
\label{bound2}
\end{equation}
\begin{equation}
L=m-\sum |R_{ijk}|.
\label{bound3}
\end{equation}
If the bounds satisfy the condition that $L \leq c \leq U$, the surface may have intersection with the cube.
Otherwise, the cube should be ruled out.

\subsubsection{Approximation by trilinear polynomial}
In each left cube, $\phi(x,y,z)$ is approximated by polynomial $\tilde P(x,y,z)$ as shown in eq \eqref{tensor1}. Then we divide each left cube into 8 smaller child cubes, and express $\tilde P(x,y,z)$ in each child cube by Legendre polynomials as follows
\begin{equation}
\tilde P(x,y,z) = A' \times_1 \vec L(\frac{2x-(a+b)}{b-a})\times_2 \vec L(\frac{2y-(c+d)}{d-c})\times_3 \vec L(\frac{2z-(e+f)}{f-e}),
\label{childtensor}
\end{equation}
where $x \in [a_1,b_1], y \in [c_1,d_1]$ and $z \in [e_1,f_1]$. Here, $[a_1,b_1] \times [c_1,d_1]\times [e_1,f_1]$ is the range of the child cube and $A'$ is the coefficient tensor of the Legendre polynomials in the child cube. When the child cubes become smaller, the coefficients of the higher order Legendre polynomials in the coefficient tensor go to zero.
The division process is repeated until the coefficients of the higher order Legendre polynomials are close to zero enough to be neglected in all the left cubes.

After above division and estimation process, in each left cube, we approximate the surface $\tilde P(x,y,z) = c$ by the following trilinear interpolation.
Supposing the range of the left cube is $[-1,1]\times[-1,1]\times[-1,1]$, the trilinear interpolation can be written in terms of the vertex values:
\begin{equation}
\begin{aligned}
g(x,y,z)= & \frac{1}{8}[\tilde P(-1,-1,-1)(1-x)(1-y)(1-z)+\tilde P(-1,-1,1)(1-x)(1-y)(1+z)\\
&+\tilde P(-1,1,-1)(1-x)(1+y)(1-z)+\tilde P(-1,1,1)(1-x)(1+y)(1+z)\\
&+\tilde P(1,-1,-1)(1+x)(1-y)(1-z)+\tilde P(1,-1,1)(1+x)(1-y)(1+z)\\
&+\tilde P(1,1,-1)(1+x)(1+y)(1-z)+\tilde P(1,1,1)(1+x)(1+y)(1+z)].
\end{aligned}
\label{linear1}
\end{equation}

\subsection{Triangulating the trilinear surface}

In this subsection, we introduce our method of triangulating the piecewise trilinear surface in cubes with different sizes. Without loss of generality, suppose in cube $[-1,1]\times[-1,1]\times[-1,1]$, the trilinear surface is $g(x,y,z) = c$.

\begin{figure}[htbp]
  \centering
  \subfigure[]{
  \label{triangulation:a}
  \includegraphics[scale = 0.7]{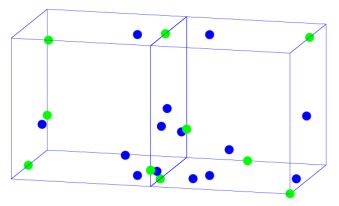}}
  \subfigure[]{
  \label{triangulation:b}
  \includegraphics[scale = 0.7]{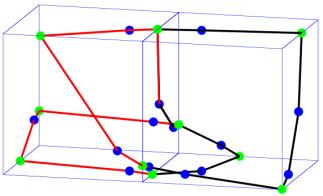}}
  \subfigure[]{
  \label{triangulation:c}
  \includegraphics[scale = 0.7]{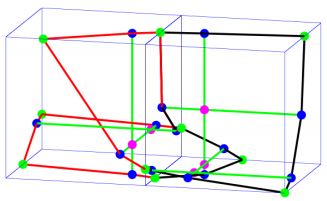}}
    \subfigure[]{
  \label{triangulation:d}
  \includegraphics[scale = 0.7]{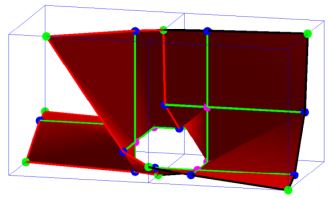}}
  \caption{Method of triangulating the trilinear surface. (a) Step 1, compute the intersection points and extreme points on the faces of the cubes. The green points are intersection points between the surface and the edge of the cube and the blue points are extreme points. (b) Step 2, connect the green intersection points and the blue extreme points by surface curves on the faces of the cubes to form red loop in the left cube and black loop in the right cube.  (c) The green lines are fold curves and the magenta points are critical points. (d) Step 3, the surface patches enclosed by the red and black loops are divided into single valued pieces by the fold curves. }
  \label{triangulation}
\end{figure}

This method contains three steps, which is shown in figure \ref{triangulation}. This figure shows the triangulation process in two neighbored cubes. Firstly, the intersection points between $g(x,y,z)=c$ and the edges of a cube are computed. They are defined as
\begin{equation}\label{intersection}
    \left \{
    \begin{aligned}
    g(x,y,z) = c \\
    \alpha = a \\
    \beta = b
    \end{aligned}
    \right.
    , \alpha, \beta \in \{x,y,z\} ,\alpha \neq \beta,  a, b \in \{1,-1\}.
  \end{equation}
The extreme points on the faces of cubes  are also computed. They are defined as
\begin{equation}\label{curve}
    \left \{
    \begin{aligned}
    g(x,y,z) = c \\
    \alpha = a \\
    \frac{\partial g(x,y,z)}{\partial \alpha} = 0
    \end{aligned}
    \right.
    , \alpha \in \{x,y,z\} , a \in \{1,-1\}.
  \end{equation}
Secondly, the intersection points and extreme points defined by eq \eqref{intersection} and eq \eqref{curve} are connected by surface curves on the faces of cube and form closed loops.
Since the surface curves on the faces of cube are simple hyperbola and the expression of the curves are explicit, it is easy to determine which two points are neighbored in the same branch of the hyperbola.
To ensure the continuity, the points belongs to the neighbored cubes and also in the current cube should be considered as well.
The surface patches enclosed by these loops may contain holes and tunnels. In the third step, the surface patches are divided into single valued pieces along $x, y, z$ directions by fold curves. Here the fold curves are defined as
\begin{equation}
\{g \left(x,y,z\right)= c, \frac{\partial g(x,y,z)}{\partial \alpha}=0\},\alpha \in \{x,y,z\}.
\label{fold1}
\end{equation}
Generally, the fold curves are not straight lines (See figure 3 in ref. \cite{TMS12}). But for the trilinear surface, the fold curves are straight line segments whose ends are extreme points.
And the fold curves along different directions may have intersections, they are critical points satisfying
\begin{equation}
    \left\{
    \begin{aligned}
    g(x,y,z) = c \\
    \frac{\partial g(x,y,z)}{\partial \alpha} = 0 \\
    \frac{\partial g(x,y,z)}{\partial \beta} = 0
    \end{aligned}
    \right.
    , \alpha, \beta \in \{x,y,z\}, \alpha \neq \beta.
    \label{critical}
\end{equation}
\begin{figure}[htbp]
  \centering
  \includegraphics[scale=0.7]{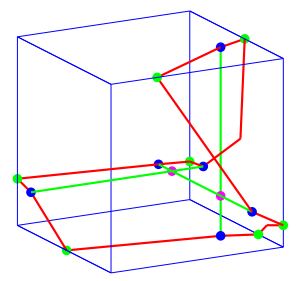}\\
  \caption{An example of connecting surface extreme points along the fold curves. The green points are intersection points and the blue points are extreme points. The green lines are fold curves. The red curve forms a close loop on the trilinear surface. The surface patches enclosed by the red loop is not single valued along $x, y, z$ directions and it is divided into four six single valued pieces along $x,y,z$ directions by the fold curves. }\label{foldcurve}
\end{figure}
Figure \ref{foldcurve} shows an example of subdividing a surface patch into single valued pieces along $x,y,z$ directions by fold curves.
The trilinear surface defined in eq \eqref{linear1} is folded at the fold curves. Cutting the trilinear surface along these fold curves ensures the resulted pieces are single valued on $x,y,z$ directions.
Subdividing the loops along fold curves helps avoid incorrect connections during triangulation and helps find missed small surface structures, such as tunnels and holes, because these structures also fold at these curves.
After the third step, each single valued piece is homomorphic to a two-dimensional polygon, and can be triangulated by standard method, such as ear clipping algorithm\cite{98triangulation,meisters1975polygons}.

\section{Experimental Results} \label{sec2}

\subsection{Efficiency and Robustness}
\begin{table}
\footnotesize
  \centering
  \begin{tabular}{c c c}
  \hline
  Molecule & & \\
  (name or PDB code) & Number of Atoms & Description \\
  \hline
  GLY & 7 & a single glycine residue \\
  ADP & 39 & ADP molecule \\
  2LWC & 75 & Met-enkephalin in DPMC SUV \\
  FAS2 & 906 & fasciculin2, a peptidic inhibitor of AChE \\
  AChE monomer & 8280 & mouse acetylcholinesterase monomer \\
  AChE tetramer & 36638 & the structure of AChE tetramer, taken from ref 29 \\
  30S ribosome & 88431 & 30S ribosome, the PDB code is 1FJF \\
  70S ribosome & 165337 & obtained from 70S\_ribosome3.7A\_model140.pdb.gz on \\
   & & http://rna.ucsc.edu/rnacenter/ribosome\_downloads.html \\
  3K1Q & 203135 & PDB code, a backbone model of an aquareovirus virion \\
  2X9XX & 510727 & a complex structure of the 70S ribosome bound to release  \\
   & & factor 2 and a substrate analog, which has 4 split PDB entries: \\
   & &  2X9R, 2X9S, 2X9T, and 2X9U \\
   1K4R & 1082160 & PDB code, the envelope protein of the dengue virus \\
  \hline
  \end{tabular}
  \caption{Description of Molecules in the PQR Benchmark}\label{pqrbench}
\end{table}
Because MSMS is the most widely used efficient software for molecular surface triangulation, in this section, the performance of TMSmesh 2.0 is compared with those of MSMS and the old version of TMSmesh.
A set of biomolecules with different sizes is chosen as a test benchmark (see Table \ref{pqrbench}) which was used in our previous work \cite{TMS11,TMS12} and can be downloaded from \emph{http://lsec.cc.ac.cn/$\sim$lubz/Download/PQR\_benchmark.tar }. The meshing softwares are run on molecular PQR files (PDB + atomic charges and radii information).
To make a reasonable comparison with MSMS, appropriate parameters, such as the error tolerance between Gaussian surface and approximated piecewise trilinear surface, are chosen for TMSmesh to achieve the surface vertex densities $1/\AA^2$ and $2/\AA^2$ used in MSMS mesh generation.
The probe radius in MSMS is set to be $1.4\AA$.
All computations run on a computer with  Intel$^\circledR$ Xeon$^\circledR$ CPU E5-4650 v2 2.4GHz and 126GB memory under 64bit Linux system.

\begin{table}
\footnotesize
  \centering
  \begin{tabular}{c c c c c c c c}
  \hline
  \multirow{2}*{Molecule} & \multirow{2}*{Natoms} & \multicolumn{3}{c}{Number of vertices} & \multicolumn{3}{c}{CPU time} \\
  \cmidrule(lr){3-5} \cmidrule(lr){6-8}
  & & TMSmesh & TMSmesh 2.0 & MSMS & TMSmesh & TMSmesh 2.0 & MSMS \\
  \hline
  FAS2 & 906 & 5170 & 6849 & 5258 & 6.4 & 0.36 & 0.13 \\
   & & 8309& 8579 & 7888 & 8 & 0.43 &0.18 \\
   AChE monomer & 8280 & 24556 & 45711 & 34819 & 52 & 1.79 & 0.72 \\
   & & 39289 & 63836 & 51784 & 60 & 2.05 & 0.96 \\
   AChE tetramer & 36638 & 95433 & 163736 & 132803 & 224 & 5.90 & 4.99 \\
   & & 152035 & 220089 & 192545 & 260 & 6.91 & 5.94 \\
   30S ribosome & 88431 & 274297 & 489325 & 353272 & 721 & 14.89 & 13.21 \\
   & & 439020 & 631448 & 520986 & 1120 & 17.59 & 15.43 \\
   70S ribosome & 165337 & 698055 & 869930 & 845550 & 1218 & 24.12 & 36.44 \\
   & & 1111399 & 1160622 & Fail & 1361 & 30.34 & Fail \\
   3K1Q & 203135 & 509390 & 678915 & 666517 & 1440 & 26.92 & 36.85 \\
   & & 812774 & 975334 & 984234 & 1728 & 30.72 & 40.48 \\
   2X9XX & 510727 & 1585434 & 2132433 & Fail & 4809 & 68.64 & Fail \\
   & & 2521233 & 2933346 & Fail & 5762 & 84.71 & Fail \\
   1K4R & 1082160 & 3325975 & 4050952 & Fail & 7296 & 141.51 & Fail \\
   & & 5298234 & 5540049 & Fail & 12905 & 178.85 & Fail \\

  \hline
  \end{tabular}

  \caption{CPU Time use for Molecular Surface Generation by TMSmesh and MSMS.}\label{CPU}
\end{table}

Table \ref{CPU} shows the CPU time cost by MSMS and TMSmesh with 1 and 2 vertex$/\AA^2$ mesh densities. In Table \ref{CPU}, TMSmesh denotes the old version in 2012\cite{TMS12} and TMSmesh 2.0 is the new version in this paper.
The discrepancies between the numbers of vertices of TMSmesh mesh and MSMS mesh are due to different definitions of molecular surface and different meshing methods used in the two programs.
The CPU time cost by TMSmesh 2.0 is much less than that cost by the old version of TMSmesh.
TMSmesh 2.0 is at least thirty times faster than the old version.
This is due to the following reasons.
Firstly, the new adaptive way of partition process to locate the surface reduces the number of surface-intersecting cubes.
We use different sizes of cubes according to the approximation accuracy of the piecewise trilinear surface in the new method instead of using same sized cubes in previous method.
Less cubes are used to precisely locate the surface.
Secondly, a more efficient and much sharper bound estimator of summation of Gaussian kernels in a cube is adopted as shown in section \ref{estimate}.
Thirdly, the trilinear polynomials are used to approximate the surface to reduce computation cost greatly.
For trilinear surface, the surface points and fold curves can be computed explicitly, and the fold curves are explicit straight lines, which make the tracing process more easily.

\begin{figure}[htbp]
  \centering
  \includegraphics[scale = 0.7]{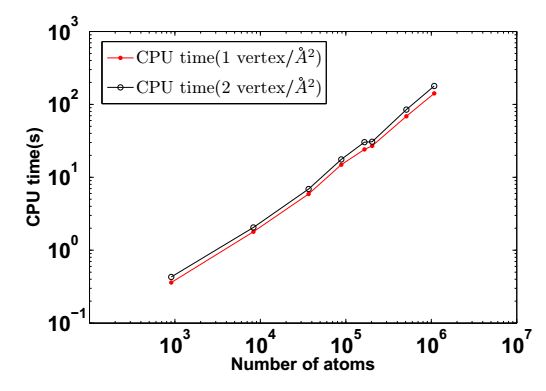}\\
  \caption{Computational performance of TMSmesh 2.0.}\label{CPUTMS}
\end{figure}
For the small molecules, the CPU time cost by MSMS is less than that of TMSmesh 2.0.
But for the large molecules, MSMS requires more time than TMSmesh 2.0.
This is because that the computational complexity of MSMS is $O[Nlog(N)]$, where $N$ is the number of atoms.
And the complexity of TMSmesh 2.0 is $O(N)$, which is shown in Figure \ref{CPUTMS}.
In TMSmesh 2.0, as the exponential kernels $e^{-||\vec x-\vec x_i||^2}$  in Gaussian surface decay very fast when the distances $||\vec x-\vec x_i||$ goes to large, all the calculations are done locally, and no global information is needed in the whole process.
TMSmesh 2.0 can successfully generate surface mesh for the biomolecules consisting of more than one million atoms, such as the dengue virus 1K4R.
Because the virus structure is among the largest ones in the Protein Data Bank, together with consideration of good algorithm stability, TMSmesh 2.0 is capable of handling the biomolecules with arbitrary sizes.

\subsection{Manifoldness}
We study the manifoldnesses of the meshes generated by TMSmesh 2.0 and MSMS. The generated surface meshes should be manifold. A non-manifold mesh can cause numerical problems in boundary element method and finite element method simulations of biomolecules.
And, non-manifold surface can not be directly used to generate the corresponding volume mesh due to its non-manifold errors, such as intersections of triangles.
The previous TMSmesh has been shown to be able to guarantee manifold mesh generation.\cite{TMS12}
Here, we check whether the meshes produced by TMSmesh 2.0 and MSMS are manifolds.
A manifold mesh for a closed molecular surface should satisfy the following three necessary conditions.\cite{TMS12}
\begin{enumerate}[(a)]
  \item Each edge should be shared and only be shared by two faces of the mesh.
  \item Each vertex should have and only have one neighborhood node loop.
  \item The mesh has no intersecting face pairs.
\end{enumerate}

Table \ref{manifold} shows the number of non-manifold defects and number of intersecting triangle pairs in the meshes produced by TMSmesh 2.0 and MSMS.
Here, the number of non-manifold defects is the number of vertices whose neighborhood does not satisfy aforementioned necessary conditions (a) and (b) for a manifold mesh.
The meshes produced by TMSmesh 2.0 all satisfy the three necessary conditions for a manifold mesh.
However, the meshes of large biomolecules generated by MSMS are not manifold.
\begin{table}
\footnotesize
  \centering
  \begin{tabular}{c c c c c c c c}
  \hline
  \multirow{2}*{Molecule} & \multirow{2}*{Natoms} & \multicolumn{2}{c}{Number of non-manifold defects} & \multicolumn{2}{c}{Number of intersecting triangle pairs} \\
  \cmidrule(lr){3-4} \cmidrule(lr){5-6}
  & & TMSmesh 2.0 & MSMS & TMSmesh 2.0 & MSMS \\
  \hline
  FAS2 & 906 & 0 & 0 & 0 & 0\\
   & & 0 & 0 & 0 & 0\\
   AChE monomer & 8280 &  0 & 0 & 0 & 220\\
   & & 0 & 0 & 0 & 265\\
   AChE tetramer & 36638 & 0 & 3 & 0 & 499\\
   & & 0 & 50 & 0 & 662\\
   30S ribosome & 88431 & 0 & 2 & 0 & 1583 \\
   & & 0 & 4 & 0 & 2504 \\
   70S ribosome & 165337 & 0 & 11 & 0 & 5235\\
   & & 0 & Fail & 0 & Fail \\
   3K1Q & 203135 & 0 & 15 & 0 & 893\\
   & & 0 & 15 & 0 & 1890 \\
   2X9XX & 510727 & 0 & Fail & 0 & Fail \\
   & & 0 & Fail & 0 & Fail \\
   1K4R & 1082160 & 0 & Fail & 0 & Fail \\
   & &  0 & Fail & 0 & Fail \\

  \hline
  \end{tabular}

  \caption{Number of non-manifold errors in meshes produced by TMSmesh 2.0 and MSMS.}\label{manifold}
\end{table}

\begin{figure}[htbp]
  \centering
  \includegraphics[scale = 1.0]{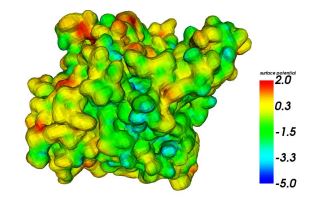}
  \caption{Electrostatic potential surface of AChE calculated by AFMPB.}\label{energy}
\end{figure}
\subsection{Boundary Element Method Simulation}
The surface mesh generated by TMSmesh 2.0 can be applied not only to molecular visualization and analysis of surface area, topology and volume in computational structure biology and structural bioinformatics, but also to boundary element method simulations.
We test the meshes in boundary element calculations of the Poission-Boltzmann electrostatics.
The BEM software used is a publicly available PB solver, AFMPB\cite{AFMPB15}.
As a representative molecular system, we choose the structure AChE monomer (see Table \ref{pqrbench}).
The surface mesh is generated by TMSmesh 2.0 and contains 87044 nodes.
Figure \ref{energy} shows the computed electrostatic potentials mapped on the molecular surface.

\subsection{Convergence}
\begin{figure}

  \includegraphics[scale = 0.55]{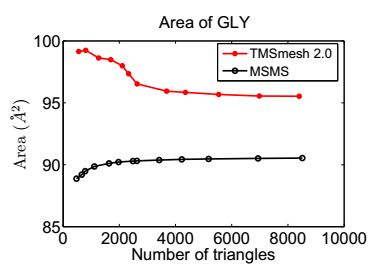}
  \includegraphics[scale = 0.55]{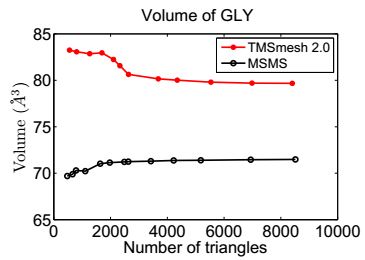}
  \includegraphics[scale = 0.55]{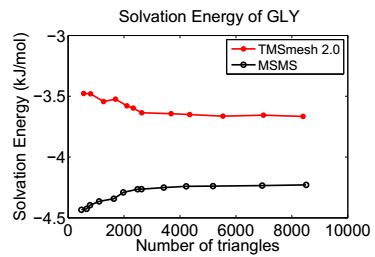}\\
  \includegraphics[scale = 0.55]{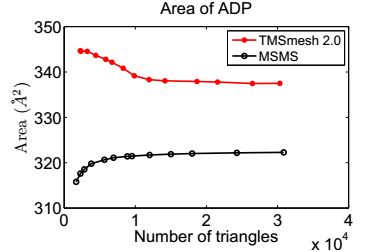}
  \includegraphics[scale = 0.55]{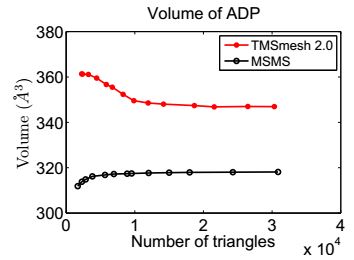}
  \includegraphics[scale = 0.55]{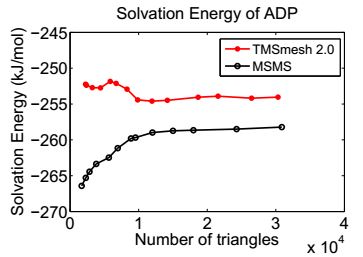}\\
  \includegraphics[scale = 0.55]{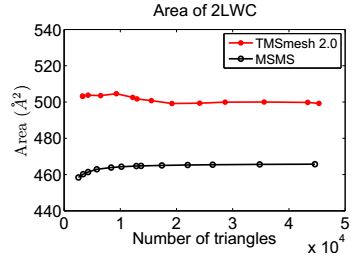}
  \includegraphics[scale = 0.55]{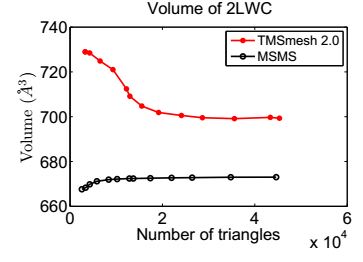}
  \includegraphics[scale = 0.55]{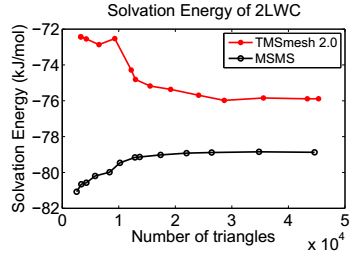}\\
  \caption{Area (left column), volume (middle column) and solvation energy (right column) for GLY (first row), ADP (second row)
and 2LWC (third row).}\label{area}
\end{figure}

Figure \ref{area} shows the solvation energies by AFMPB as well as the surface areas and molecular volumes computed from the meshes of three small molecules, GLY, ADP and 2LWC (see Table \ref{pqrbench}) using different mesh densities.
The results show that the meshes produced by TMSmesh 2.0 lead to convergent and reasonable results for energy, area and volume when the mesh density increasing.
However, the results computed by MSMS converge a little more smoothly than those of TMSmesh 2.0 when the number of triangles are not large.
This is because that we use the trilinear polynomial to approximate the Gaussian kernel function.
Less triangles lead to lower precisions of the approximation, which causes more uncertainties.
The disparities between the limits when number of triangles goes to large are due to the different molecular surface definitions used by TMSmesh and MSMS.

\subsection{Volume Mesh Generation Conforming Surface Mesh}

The surface mesh generated by TMSmesh 2.0 can be directly used to generate corresponding surface conforming volume mesh.
And the volume mesh generated by this method can be applied to the finie element method simulation directly.
\begin{figure}
  \centering
  \includegraphics[scale = 0.8]{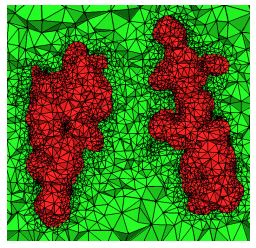}\\
  \caption{A cross section of the volume mesh for VDAC. The surface mesh is generated by TMSmesh 2.0 and the volume mesh is generated by TetGen.}\label{2JK4}
\end{figure}
Figure \ref{2JK4} shows a cross section of the volume mesh for the ion channel, VDAC (PDB code: 2JK4).
The VDAC serves an essential role in the transport of metabolites and electrolytes between the cell matrix and mitochondria.\cite{vdac}
For this example, the molecular surface mesh is generated by TMSmesh 2.0 and the corresponding volume mesh is generated by TetGen\cite{tetgen2}.
The channel pore is clearly represented in the mesh and the detailed topology is correctly preserved, which is important for ion channel simulations.
In addition, from the cross section we can see that the surface mesh is dense at the rugged parts and sparse at the smooth parts.


\section{ Conclusion}\label{sec3}

We have described a new algorithm in TMSmesh 2.0 for triangulating the Gaussian molecular surface.
In TMSmesh 2.0, an adaptive surface partition is developed using a new method to estimate the upper and lower bounds of surface function in a cell. In each located cube, a trilinear polynomial is used to approximate the Gaussian surface within controllable precision.
The fold curves are used to divide the trilinear surface in each cube into single valued pieces to guarantee a manifold mesh generation.
Compared with the old version, TMSmesh 2.0 is more than thirty times faster.
TMSmesh 2.0 is shown to be a robust and efficient software to mesh the Gaussian molecular surface.
The meshes generated by TMSmesh 2.0 are manifold without intersections.
And the mesh can be directly used in boundary element type of simulations and volume mesh generations.

\section{ Acknowledgements}
Tiantian Liu and Benzhuo Lu are supported by the State Key Laboratory of Scientific/Engineering
Computing, National Center for Mathematics and Interdisciplinary
Sciences, Science Challenge Program (SCP) and the China NSF (NSFC 91530102, NSFC 21573274).
Minxin Chen is supported by China NSF (NSFC11301368) and the NSF of Jiangsu Province (BK20130278).

\bibliographystyle{unsrt}
\bibliography{reference}
\clearpage

\end{document}